\documentclass[twoside,reqno,11pt]{fcaa-var} %
\usepackage{graphicx}
\usepackage{epsfig}
\usepackage{amsthm}
\usepackage{amsmath}
\usepackage{latexsym}
\usepackage{amsfonts}
\usepackage{amssymb}

\newtheoremstyle{theorem}
  {15pt}          
  {15pt}  
  {\sl}  
  {\parindent}
  {\sc}  
  {. }   
  { }    
  {}     
\theoremstyle{theorem}
\newtheorem{lemma}{Lemma}[section]

\newtheoremstyle{defi}
  {15pt}          
  {15pt}  
  {\rm}  
  {\parindent}     
  {\sc}  
  {. }    
  { }    
  {}     
\theoremstyle{defi}

 \def\proofname{\indent\rm P\ r\ o\ o\ f.}
 \def\proofend{\hfill$\Box$}
 \def\qed{\hfill$\Box$} 

 \usepackage{hyperref} 

 \def\theequation{\arabic{section}.\arabic{equation}}

 \def\themycopyrightfootnote{\vspace*{3pt}
 Submitted to Fract. Calc. Appl. Anal., 
 \par \url{https://www.degruyter.com/view/j/fca}
  }

  \newtheorem{result}[lemma]{Standard Result}

 \title[The Cauchy derivative as a principal value]{The generalised Cauchy derivative as a principal value of the Gr\"unwald-Letnikov fractional derivative  for divergent expansions}
 \author[\normalsize A. Pallavi Sudhir]{\normalsize Abhimanyu Pallavi Sudhir $^1$}

\begin{document}

 \vbox to 2.5cm { \vfill }


 \bigskip \medskip

 \begin{abstract}

It has recently been proven that the generalised Cauchy fractional derivative (also known as the Riemann-Liouville fractional derivative) is equal to the Gr\"unwald-Letnikov derivative. However, we observe that there are ``Gr\"unwald non-differentiable" functions for which the latter derivative is not convergent, while the Riemann-Liouville derivative is. In this paper, we show that the Riemann-Liovuille derivative can be considered a ``principal value" of the Gr\"unwald-Letnikov derivative, requiring specific relative rates of approach for the limits $\Delta x \to 0$ and $N \to \infty$ (where $N$ is the upper limit of the infinite summation) in the Gr\"unwald derivative -- i.e. instead of varying $\Delta x$ and $N$ independently, the two must satisfy a relation and be varied as a single limit. We proceed to calculate this relation for several functions and orders, finding that several possible such relations are possible for a given function, also placing requirements on the ``handedness" of the fractional derivative. It is further shown that for functions with a Taylor expansion, the relation $\Delta x=x/N$ always produces the correct principal value.
 \medskip


 \smallskip


 \end{abstract}

 \maketitle

 \vspace*{-16pt}



 \section{Introduction}\label{sec:1}

\setcounter{section}{1}
\setcounter{equation}{0}\setcounter{theorem}{0}

The Gr\"unwald-Letnikov fractional derivative is a classical fractional derivative that arises from the generalisation of the limit form of the $n$th order derivative, by replacing the finite summation with an infinite summation, and can be written as (among other equivalent definitions):
\begin{equation}\label{gl-def}
D^Rf(x)=\lim\limits_{h\to0}h^{-R}\sum\limits_{k=0}^\infty \left(\begin{array}{*{20}{c}} {R} \\ k \end{array}\right)(-1)^kf(x-kh)
\end{equation}

The equivalence of this derivative to the standard Riemann-Liouville derivative has been a matter of considerable mathematical debate and discussion in the past. It was conjectured in \cite{diaz1974} that these two classical fractional derivatives are, in fact, equal. This was discussed, e.g. in \cite{kiry1996}, before eventually being proven in \cite{orti2004}.

However, it is easy to find examples of functions for which the Gr\"unwald-Letnikov derivative \emph{isn't} well-defined (or at least not uniquely defined) -- for example, finite-degree non-constant polynomials. For the simplest example, consider the half-derivative, $D^{1/2}$ of $f(x)=x$. Here, the Gr\"unwald derivative is calculated as:
\begin{equation}\label{halfdiverge}
\begin{gathered}
  {D^{1/2}}x = \mathop {\lim }\limits_{h \to 0} {h^{ - 1/2}}\sum\limits_{k = 0}^\infty  {\left( {\begin{array}{*{20}{c}}
  {1/2} \\ 
  k 
\end{array}} \right)} {( - 1)^k}(x - kh) \\ 
   = {h^{ - 1/2}}\left[ {x\sum\limits_{k = 0}^\infty  {\left( {\begin{array}{*{20}{c}}
  {1/2} \\ 
  k 
\end{array}} \right)} {{( - 1)}^k} - h\sum\limits_{k = 0}^\infty  {\left( {\begin{array}{*{20}{c}}
  {1/2} \\ 
  k 
\end{array}} \right)k} {{( - 1)}^k}} \right] \\ 
   = \frac{{\sqrt h }}{2}\sum\limits_{k = 0}^\infty  {\left( {\begin{array}{*{20}{c}}
  { - 1/2} \\ 
  k 
\end{array}} \right)} {( - 1)^k} \\ 
\end{gathered} 
\end{equation}

Where in the last step we used the fact that $\sum_{k = 0}^\infty  {\left( {\begin{array}{*{20}{c}}
  {1/2} \\ 
  k 
\end{array}} \right)k} {( - 1)^k} = 0$ for $k=0$, then transformed the limits of the summation. The resulting sum in Eq.~\eqref{halfdiverge} is clearly outside the interval of convergence of the binomial series (one may informally regard it as the binomial series for $(1-1)^{-1/2}$). Meanwhile, the Riemann-Liouville derivative of a power function is simply:
\begin{equation}\label{rl-def}
D^{R}x^m=\frac{\Gamma(m+1)}{\Gamma(m-R+1)}x^{m-R}
\end{equation}

Which for $m=1$, $R=1/2$ evaluates to $\frac{2}{\sqrt\pi}\sqrt{x}$. This breaking of the proven equivalence may seem problematic, but it is important to note that the divergence was not of the Gr\"unwald-Letnikov derivative per se, but only of the infinite sum. This divergent value is then multiplied by $h$, that approaches zero, yielding an undefined answer.

Noting that an infinite sum is simply the limit of the partial sum, i.e. $\sum_{k=0}^\infty=\lim_{N\to\infty}\sum_{k=0}^N$ We are thus led to ask if we vary $N$ and $h$ together in a specific way, we can recover the Riemann-Liouville fractional derivative as a ``principal value" of the Gr\"unwald derivative.

   \vspace*{12pt} 
 \subsection{Standard results and lemmas}\label{subsec:1.1}

This section contains standard results pertaining to hypergeometric and related special functions, and generalised Newtonian binomial coefficients and summations involving them, included here for reference in later proofs. They can be found in most elementary textbooks on the subject, e.g. \cite{gui2018}.

\begin{result}\label{lm1}
A ``generalised Pascal rule" for fractional $R$:
\begin{equation*}
\sum\limits_{k = 0}^N {\left( {\begin{array}{*{20}{c}}
  {R} \\ 
  k 
\end{array}} \right){{( - 1)}^k}}  = {( - 1)^N}\left( {\begin{array}{*{20}{c}}
  {R - 1} \\ 
  N 
\end{array}} \right)\end{equation*}
 \end{result}

\begin{result}\label{lm2}
\begin{equation*}
\sum\limits_{j = 0}^m {\left( {\begin{array}{*{20}{c}}
  m \\ 
  j 
\end{array}} \right)\frac{{{{\left( { - q} \right)}^j}}}{{R - j}}}  = {\frac{{_2F_1}\left( { - m, - R; \, 1 - R; \, q} \right)}{R}}
\end{equation*}
\end{result}

\begin{result}\label{lm3}
\begin{equation*}
\sum\limits_{j = 0}^m {\left( {\begin{array}{*{20}{c}}
  m \\ 
  j 
\end{array}} \right)\frac{{{{\left( { - 1} \right)}^j}}}{{R - j}}}  = {{R^{ - 1}}{{\left( {\begin{array}{*{20}{c}}
  {m - R} \\ 
  m 
\end{array}} \right)}^{ - 1}}}
\end{equation*}
\end{result}

\begin{result}\label{lm4}
For all $R$,
\begin{equation*}
    \Gamma(1+R)\Gamma(1-R)=\frac{\pi R}{\sin(\pi R)}
\end{equation*}
\end{result}

\begin{result}\label{lm5}
As $N\to\infty$, 
\begin{equation*}
{\left( { - 1} \right)^N}\left( {\begin{array}{*{20}{c}}
  R \\ 
  N 
\end{array}} \right) \sim -\frac{\sin(\pi R)}{\pi}\,\Gamma(R+1)\,{N^{ - R - 1}}
\end{equation*}
\end{result}

\begin{result}\label{lm6}
As $m\to\infty$, 
\begin{equation*}
\left( {\begin{array}{*{20}{c}}
  m \\ 
  R 
\end{array}} \right) \sim  - \frac{\sin(\pi R)}{\pi }\,\Gamma ( - R)\,{m^R}
\end{equation*}
\end{result}

\begin{lemma}\label{lm7}
As $N\to\infty$, for integer $j$
\begin{equation*}
\sum\limits_{k = 0}^N {\left( {\begin{array}{*{20}{c}}
  R \\ 
  k 
\end{array}} \right){{\left( { - 1} \right)}^k}{k^j}}  \sim {\left( { - 1} \right)^N}\frac{R}{{R - j}}\left( {\begin{array}{*{20}{c}}
  {R - 1} \\ 
  N 
\end{array}} \right){N^j}
\end{equation*}
\end{lemma}

\proof 
 The actual value of the summation in Lemma~$\eqref{lm7}$ in general is a $j$th order polynomial in $N$, but the lower-order terms vanish for large $N$. This can be proven fairly easily via induction in $j$ -- observing that the case $j=0$ is simply Standard Result~\ref{lm1},
 \begin{equation*}
 \begin{gathered}
  \sum\limits_{k = 0}^N {{{\left( { - 1} \right)}^k}\left( {\begin{array}{*{20}{c}}
  R \\ 
  k 
\end{array}} \right){k^{j + 1}}}  = \sum\limits_{k = 0}^N {{{\left( { - 1} \right)}^k}\frac{{\Gamma (R + 1)}}{{\Gamma (R - k + 1)\Gamma (k)}}{k^j}}  \\ 
   =  - R\sum\limits_{k = 0}^N {{{\left( { - 1} \right)}^{k - 1}}\left( {\begin{array}{*{20}{c}}
  {R - 1} \\ 
  {k - 1} 
\end{array}} \right){k^j}}  \\ 
   =  - R\sum\limits_{k =  - 1}^{N - 1} {{{\left( { - 1} \right)}^k}\left( {\begin{array}{*{20}{c}}
  {R - 1} \\ 
  k 
\end{array}} \right){{\left( {k + 1} \right)}^j}}  \\ 
   =  - R\sum\limits_{k = 0}^{N - 1} {{{\left( { - 1} \right)}^k}\left( {\begin{array}{*{20}{c}}
  {R - 1} \\ 
  k 
\end{array}} \right){k^j}}  - Rj\sum\limits_{k =  - 1}^N {{{\left( { - 1} \right)}^k}\left( {\begin{array}{*{20}{c}}
  {R - 1} \\ 
  k 
\end{array}} \right){k^{j - 1}}}  - ... \\ 
   =  - R\left[ {{{\left( { - 1} \right)}^{N - 1}}\frac{{R - 1}}{{R - 1 - j}}\left( {\begin{array}{*{20}{c}}
  {R - 2} \\ 
  {N - 1} 
\end{array}} \right){{\left( {N - 1} \right)}^j}} \right] - \left( {...} \right){\left( {N - 1} \right)^{j - 1}} - ... \\ 
   = {\left( { - 1} \right)^N}\frac{{R(R - 1)}}{{R - 1 - j}}\frac{{\Gamma (R - 1)}}{{\Gamma (N)\Gamma (R - N)}}\left( {{N^j} - j{N^{j - 1}} + ...} \right) - ... \\ 
   = {\left( { - 1} \right)^N}\frac{R}{{R - 1 - j}}\frac{{\Gamma (R)}}{{\Gamma (N + 1)\Gamma (R - N)}}N\left( {{N^j} - j{N^{j - 1}} + ...} \right) - ... \\ 
  \sim{\left( { - 1} \right)^N}\frac{R}{{R - \left( {j + 1} \right)}}\left( {\begin{array}{*{20}{c}}
  {R - 1} \\ 
  N 
\end{array}} \right){N^{j + 1}} \\ 
\end{gathered} 
\end{equation*}
Completing our proof.
 \proofend 

\section{The half-derivative of a linear function}\label{sec:2}

\setcounter{section}{2}
\setcounter{equation}{0}\setcounter{theorem}{0}

We approach our goal of obtaining the relative rates of $h\to0$ and $N\to\infty$ for the fractional derivative of a general power function by first considering the half-derivative of a simple linear function $f(x)=x$ (the result will clearly also apply to a general $f(x)=ax+b$, since the half-derivative of a constant converges). Substituting into Eq.~\eqref{gl-def} and expanding:
\begin{equation*}
{D^{1/2}}x = \lim\limits_{h\to0,\,N\to\infty}\frac{x}{{\sqrt h }}\sum\limits_{k = 0}^N {\left( {\begin{array}{*{20}{c}}
  {1/2} \\ 
  k 
\end{array}} \right){{\left( { - 1} \right)}^k}}  - \sqrt h \sum\limits_{k = 0}^N {k\left( {\begin{array}{*{20}{c}}
  {1/2} \\ 
  k 
\end{array}} \right){{\left( { - 1} \right)}^k}} 
\end{equation*}

We apply Standard Result~\ref{lm1} and Lemma~\ref{lm7} respectively to the two sums and simplify to obtain the closed form expression:
\begin{equation}\label{halfform}
{D^{1/2}}x = \frac{x}{{\sqrt h }}{\left( { - 1} \right)^N}\left( {\begin{array}{*{20}{c}}
  { - 1/2} \\ 
  N 
\end{array}} \right) + \frac{{\sqrt h }}{2}{\left( { - 1} \right)^N}\left( {\begin{array}{*{20}{c}}
  { - 3/2} \\ 
  N 
\end{array}} \right)
\end{equation}

Recall from Standard Result~\ref{lm5} that for large $N$:
\begin{equation*}
\left(-1\right)^N\left( {\begin{array}{*{20}{c}}
  { - 1/2} \\ 
  N 
\end{array}} \right)\sim\frac{1}{{\sqrt \pi  }}{N^{ - 1/2}}    
\end{equation*}
\begin{equation*}
\left(-1\right)^N\left( {\begin{array}{*{20}{c}}
  { - 3/2} \\ 
  N 
\end{array}} \right)\sim \frac{2}{{\sqrt \pi  }}{N^{1/2}}    
\end{equation*}

Thus, we may simplify \eqref{halfform} in the limit as $N\to\infty$ as:
\begin{equation}\label{halfformbign}
{D^{1/2}}x = \frac{x}{{\sqrt h }}\frac{1}{{\sqrt \pi  }}{N^{ - 1/2}} + \frac{{\sqrt h }}{2}\frac{2}{{\sqrt \pi  }}{N^{1/2}}
\end{equation}

Which must be equal in the limit to the Riemann-Liouville derivative, which can be calculated from Eq.~\eqref{rl-def}:
\begin{equation}\label{halfcondition}
\frac{2}{{\sqrt \pi  }}\sqrt x  = \frac{x}{{\sqrt h }}\frac{1}{{\sqrt \pi  }}{N^{ - 1/2}} + \frac{\sqrt{h}}{{\sqrt \pi  }}{N^{1/2}}
\end{equation}

This can be regarded as an equation in $h$, and its solution is:
\begin{equation}\label{halfsolution}
h=x/N
\end{equation}

In other words, the half-derivative should be evaluated with $h$ varying as $x/N$, which indeed approaches 0 as $N\to\infty$. This was an extraordinarily simple solution -- in a sense, the simplest possible combination of the relative rates of the two limits. As it turns out, the solution will not be so simple for the general Gr\"unwald-Letnikov derivative $D^R$ of any power function, however it will continue to take the form $h=qx/N$ for a constant $q$.

\section{The general fractional derivative of a linear function}\label{sec:3}

\setcounter{section}{3}
\setcounter{equation}{0}\setcounter{theorem}{0}

We will repeat the process we demonstrated in Sec~\ref{sec:2} to obtain the relation between $h$ and $N$ for a slightly more complicated derivative -- the $R$th derivative of $f(x)=x$. We start by applying Eq.~\eqref{gl-def} to the function, then applying Standard Result~\ref{lm1} and Lemma~\ref{lm7}:
\begin{equation*}
{D^R}x = {h^{ - R}}\left[ {x\sum\limits_{k = 0}^N {\left( {\begin{array}{*{20}{c}}
  R \\ 
  k 
\end{array}} \right){{\left( { - 1} \right)}^k}}  - h\sum\limits_{k = 0}^N {k\left( {\begin{array}{*{20}{c}}
  R \\ 
  k 
\end{array}} \right){{\left( { - 1} \right)}^k}} } \right]
\end{equation*}
\begin{equation}\label{Rform}
{D^R}x = {h^{ - R}}{\left( { - 1} \right)^N}\left( {\begin{array}{*{20}{c}}
  {R - 1} \\ 
  N 
\end{array}} \right)\left[ {x - \frac{R}{{R - 1}}Nh} \right]
\end{equation}

From Standard Result~\ref{lm5}, we have:
\begin{equation*}
{( - 1)^N}\left( {\begin{array}{*{20}{c}}
  {R - 1} \\ 
  N 
\end{array}} \right)\sim\frac{\sin(\pi R)}{\pi }\Gamma (R){N^{ - R}}
\end{equation*}

Equating the \eqref{Rform} to the definition of the Riemann-Liouville derivative in \eqref{rl-def}, it is determined that the principal value is achieved under the condition:
\begin{equation}\label{Rformbign}
\frac{1}{{\Gamma (2 - R)}}{x^{1 - R}} = \frac{\sin(\pi R)}{\pi }\Gamma (R){(Nh)^{ - R}}\left[ {x - \frac{R}{{R - 1}}Nh} \right]
\end{equation}

This may be algebraically manipulated into an equation in $Nh/x$:
\begin{equation}\label{Rcondition}
\frac{R}{{R - 1}}{(Nh/x)^{1 - R}} - {(Nh/x)^{ - R}} + \frac{\pi/\sin(\pi R)}{{\Gamma (2 - R)\Gamma (R)}} = 0
\end{equation}

In other words, $Nh/x$ must, indeed as we suspected, be a constant. Letting $q=Nh/x$, our condition for achieving the principal value of the Gr\"unwald-Letnikov fractional derivative of $f(x)=x$ is:
\begin{equation}\label{Rsolution}
h=qx/N
\end{equation}

Where $q$ depends on $R$ as a root of the ``characteristic polynomial" derived as a simplification of Eq.~\eqref{Rcondition} (via Standard Result~\ref{lm4}):
\begin{equation}\label{Rpolynomial}
{q^R} - Rq + (R - 1) = 0
\end{equation}

Some special cases of this result are immediately obvious:
\begin{itemize}
    \item $q=1$ is always a solution, implying the relation $h=x/N$ yields the principal value for $D^Rx$ for all $R$.
    \item Setting $R=1$ allows $q=0$, i.e. the limit in $N$ is taken before the limit in $h$. This is indeed the typical approach taken to take an integer-order derivative of $f(x)=x$.
    \item As an example, consider $R=1/3$. While $q=1$ is still a solution, of course, so is $q=-8$. This is in fact the value of $q$ required for the principal value of the \emph{reverse Gr\"unwald-Letnikov derivative}.
\end{itemize}

N.B. We're referring to Eq.~\eqref{Rpolynomial} as a ``polynomial" for convenience even though the order of the ``polynomial" is fractional. For rational $R$, the equation can indeed be transformed into a polynomial.

\vspace*{2pt} 
\section{The general fractional derivative of a power function}\label{sec:4}

\setcounter{section}{4}
\setcounter{equation}{0}\setcounter{theorem}{0}

This section considers the general fractional derivative of a power function, $D^Rx^m$  where we consider only integer $m$ to avoid having to use Newton's binomial theorem (which would introduce a new limit). We substitute the function into Eq.~\eqref{gl-def} and apply the binomial theorem:
\begin{equation*}
{D^R}{x^m} = \lim\limits_{h\to0,\,N\to\infty}{h^{ - R}}\sum\limits_{k = 0}^N  {\left( {\begin{array}{*{20}{c}}
  R \\ 
  k 
\end{array}} \right)} {( - 1)^k}\sum\limits_{j = 0}^m {\left( {\begin{array}{*{20}{c}}
  m \\ 
  j 
\end{array}} \right){x^{m - j}}{{( - kh)}^j}} 
\end{equation*}

We evaluate this sum via much of the same means as with the linear function in \ref{sec:3}. Reordering the summations,

\begin{equation*}{D^R}{x^m} = {h^{ - R}}\sum\limits_{j = 0}^m {{{( - 1)}^j}\left( {\begin{array}{*{20}{c}}
  m \\ 
  j 
\end{array}} \right){x^{m - j}}{h^j}\sum\limits_{k = 0}^N {\left( {\begin{array}{*{20}{c}}
  R \\ 
  k 
\end{array}} \right){{( - 1)}^k}{{k}^j}} }
\end{equation*}

The summation in $k$ is given by Lemma~\ref{lm7}, simplifying our expression to:
\begin{equation*}
\begin{gathered}
  {D^R}{x^m} = {h^{ - R}}\sum\limits_{j = 0}^m {{{( - 1)}^j}\left( {\begin{array}{*{20}{c}}
  m \\ 
  j 
\end{array}} \right){x^{m - j}}{h^j}{{( - 1)}^N}\frac{R}{{R - j}}\left( {\begin{array}{*{20}{c}}
  {R - 1} \\ 
  N 
\end{array}} \right){N^j}}  \\ 
   = {h^{ - R}}{( - 1)^N}\left( {\begin{array}{*{20}{c}}
  {R - 1} \\ 
  N 
\end{array}} \right){x^m}\sum\limits_{j = 0}^m {{{( - 1)}^j}\left( {\begin{array}{*{20}{c}}
  m \\ 
  j 
\end{array}} \right){{(Nh/x)}^j}\frac{R}{{R - j}}}  \\ 
  \sim\frac{\sin(\pi R)}{\pi }\Gamma (R + 1){(Nh/x)^{ - R}}{x^{m - R}}\sum\limits_{j = 0}^m {{{( - 1)}^j}\left( {\begin{array}{*{20}{c}}
  m \\ 
  j 
\end{array}} \right){{(Nh/x)}^j}\frac{1}{{R - j}}}  \\ 
\end{gathered} 
\end{equation*}

Where we used Standard Result~\ref{lm5}. Substituting $q=Nh/x$, equating to the Riemann-Liouville derivative in Eq.~\eqref{rl-def} and simplifying:
\begin{equation}\label{Rmpolynomial}
{q^{ - R}}\sum\limits_{j = 0}^m {\left( {\begin{array}{*{20}{c}}
  m \\ 
  j 
\end{array}} \right)\frac{{{{( - q)}^j}}}{{R - j}}}  = \frac{\pi}{\sin(\pi R)} \left( {\begin{array}{*{20}{c}}
  m \\ 
  R 
\end{array}} \right)
\end{equation}

Or in terms of hypergeometric functions via Lemma~\ref{lm2}:
\begin{equation}\label{Rmpolynomialhyper}
{q^{ - R}}_2{F_1}( - m, - R;1 - R;q) = \frac{\pi R}{\sin(\pi R)} \left( {\begin{array}{*{20}{c}}
  m \\ 
  R 
\end{array}} \right)
\end{equation}

Eqs.~\eqref{Rmpolynomial} and \eqref{Rmpolynomialhyper} are thus two alternative forms of the characteristic polynomial for $D^R x^m$, whose solution defines the relation between $h$ and $N$ via Eq.~\eqref{Rsolution}. 
\vspace*{6pt} 

One may verify that $q=1$ is still a solution to this characteristic equation. From Standard Result~\ref{lm3}, we have for $q=1$:
\begin{equation*}
{q^{ - R}}\sum\limits_{j = 0}^m {\left( {\begin{array}{*{20}{c}}
  m \\ 
  j 
\end{array}} \right)\frac{{{{( - 1)}^j}}}{{R - j}}}  = {{R^{-1}\left( {\begin{array}{*{20}{c}}
  {m - R} \\ 
  m 
\end{array}} \right)^{-1}}}
\end{equation*}

And indeed, by Standard Result~\ref{lm4},
\begin{equation*}
{R^{ - 1}}{\left( {\begin{array}{*{20}{c}}
  {m - R} \\ 
  m 
\end{array}} \right)^{ - 1}} = \frac{{\pi R}}{{\sin (\pi R)}}\left( {\begin{array}{*{20}{c}}
  m \\ 
  R 
\end{array}} \right)
\end{equation*}

As a result, setting $h=x/N$ returns the principal value of the Gr\"unwald-Letnikov derivative not only for the power function $D^Rx^m$, but for \emph{any polynomial function}, and thus for any function with a Taylor expansion.

\section{Conclusion}\label{sec:5}

\setcounter{section}{5}
\setcounter{equation}{0}\setcounter{theorem}{0}

We have shown that the Riemann-Liouville derivative can be recovered as a principal value of the Gr\"unwald-Letnikov derivative by setting $h=qx/N$ where $q$ is a root of the function's ``characteristic polynomial", defined for power functions via Eq.~\eqref{Rmpolynomialhyper}. In addition, it was demonstrated that the root $q=1$, i.e. $h=x/N$, works for any fractional derivative of any power function, and by extension any polynomial or function with a Taylor expansion.  

Our analysis reveals not only what rate the limit $h\to 0$ must be taken, but also \emph{which}-handed Gr\"unwald derivative to choose in order to recover the principal value. For instance, the reverse derivative is an alternative form of the Gr\"unwald derivative can be shown to be equivalent to it barring a transformation of $h\to-h$. However, we have shown that the principal value of the derivative requires specific rates and directions of approach for $h$. For example, $h=x/N$ is the only solution while evaluating $D^{1/2}x$ -- since $N$ is a positive infinity, evaluating the principal value requires $h$ to have the same sign as $x$. On the other hand, $D^{1/3}x$, as we saw in Sec~\ref{sec:3}, allows for both-handed derivatives to exist as solutions, $h=x/N$ and $h=-8x/N$. Similarly, it might be possible for $h$ to have a complex argument, depending on the characteristic polynomial. 

Future research may focus on finding the characteristic polynomial of a general polynomial function -- while we have established that $h=x/N$ is a universal solution, other rates of approach clearly must exist.




 \bigskip \smallskip

 \it

 \noindent
$^1$ Department of Mathematics \\
Imperial College of London \\
180 Queen's Gate, South Kensington Campus \\
London -- SW7 2AZ, UNITED KINGDOM  \\[4pt]
  e-mail: ap6218@imperial.ac.uk
\hfill Received: September 12, 2018 \\[12pt]

\end{document}